\documentclass[smallextended,envcountsect,natbib]{svjour3}       
\smartqed  
\usepackage{graphicx}
\usepackage{mathtools,amssymb}
\usepackage{xcolor}
\usepackage{longtable}
\begin{document}

\title{Determining the optimal piecewise constant approximation for the Nonhomogeneous Poisson Process rate of Emergency Department patient arrivals}



\titlerunning{Optimal piecewice constant approximation for NHPP rate of ED arrivals}        


\author{Alberto~De~Santis~\href{https://orcid.org/0000-0001-5175-4951}{\includegraphics[scale=0.5]{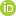}}
	\and 
Tommaso~Giovannelli~\href{https://orcid.org/0000-0002-1436-5348}{\includegraphics[scale=0.5]{orcid.png}}
	\and 
Stefano~Lucidi~\href{https://orcid.org/0000-0003-4356-7958}{\includegraphics[scale=0.5]{orcid.png}}
    \and
Mauro~Messedaglia
    \and
Massimo~Roma~\href{https://orcid.org/0000-0002-9858-3616}{\includegraphics[scale=0.5]{orcid.png}}
}

\authorrunning{A. De Santis, T. Giovannelli, S. Lucidi, M. Messedaglia, M. Roma} 


\institute{A. De Santis, T. Giovannelli, S. Lucidi, M. Roma \at
	Dipartimento di Ingegneria Informatica, Automatica e Gestionale  ``A. Ruberti'' \\
		SAPIENZA, Universit\`a di Roma \\ 
		via Ariosto, 25 -- 00185 Roma, Italy \\
		\email{desantis,giovannelli,lucidi,roma@diag.uniroma1.it}
           \and
           M. Messedaglia \at
              ACTOR Start up of SAPIENZA Universit\`a di Roma \\
              via Nizza 45, 00198 Roma, Italy.\\
              \email{mauro.messedaglia@gmail.com}
}

\date{}

\maketitle


\begin{abstract}
Modeling the arrival process to an Emergency Department (ED) is the first step of all studies dealing with the patient flow within the ED. Many of them focus on the increasing phenomenon of ED overcrowding, which is afflicting hospitals all over the world. Since Discrete Event Simulation models are often adopted with the aim to assess solutions for reducing the impact of this problem, proper nonstationary processes are taken into account to reproduce time-dependent arrivals. Accordingly, an accurate estimation of the unknown arrival rate is required to guarantee reliability of results.\\
In this work, an integer nonlinear black-box optimization problem is solved to determine the 
best piecewise constant approximation of the time-varying arrival rate function, by finding the optimal partition of the 24 hours into a suitable number of non equally spaced intervals. The black-box constraints of the optimization problem make the feasible solutions satisfy proper statistical hypotheses; these ensure the validity of the nonhomogeneous Poisson assumption about the arrival process, commonly adopted in the literature, and prevent to mix overdispersed data for model estimation. The cost function includes a fit error term for the solution accuracy and a penalty term to select an adeguate degree of regularity of the optimal solution. To show the effectiveness of this methodology, real data from one of the largest Italian hospital EDs are used. 	\\	
\keywords{Emergency Department \and Arrival process \and Non Homogeneous Poisson Process \and Black-Box Optimization}
\end{abstract}

\section{Introduction}\label{sec:intro}
Statistical modelling for describing and predicting patient arrival to Emergency Departments (EDs) represent a basic tool of each study concerning ED patient load and crowding. Indeed, all the approaches adopted to this aim require an accurate model of the patient arrival process. Of course, such a process plays a key role in tackling the widespread phenomenon of overcrowding which afflicts EDs all over the world \cite{Ahalt.2018,Bernstein.2003,Daldoul.2018,hoot.aronsky:08,Hoot.2007,Reeder.2003,Vanbrabant:20,Wang.2015,Weiss.2004,Weiss.2006}. 
The two factors that have the most significant effect on overcrowding are both external and internal. The first concerns the patient arrival process; the second regards the patient flow within the ED. Therefore, both the aspects must be accurately considered for a reliable study on ED operation.
\par
Several modelling approaches for analyzing ED patient flow have been proposed in literature (see \cite{Wiler:2011} for a survey). The main quantitative methods used are based on statistical analysis (time--series, regression) or on general analytic formulas (queuing theory). Simulation modelling (both Discrete Event and Agent Based Simulation) is currently one of the most widely used and flexible tool for studying the patient flow through an ED. In fact, it enables to perform an effective scenario analysis, aiming at determining bottlenecks (if any) and testing different ED settings. We refer to \cite{Salmon:2018} for a recent survey on simulation modelling for ED operation. A step forward is Simulation--Based Optimization methodology which combines a simulation model with a black-box optimization algorithm, aiming at determining an optimal ED settings, based on suited objective function (representing some KPIs) to be maximized or minimized \cite{Ahmed.2009,Guo.2016,Guo.2017}.
\par
Modeling methodologies are generally based on assumptions that, in some cases, may represent serious limitations when applied to complex real--world cases, such as ED operation. In particular, when dealing with ED patient arrival stochastic modelling, due to the
nonstationarity of the process, a standard assumption is the use of {\em Nonhomogeneous Poisson Process (NHPP)} \cite{Ahalt.2018,Ahmed.2009,Guo.2017,Kim.2014,Kuo.2016,Zeinali.2015}. We recall that a counting process $X(t)$ is a NHPP if {\em 1)} arrivals occur one at a time (no batch); {\em 2)} the process has independent increments; {\em 3)} increments have Poisson distribution, i.e. for each interval $[t_1, t_2]$,
	$$
	P\left( X(t_1)-X(t_2)=n\right) = e^{-m(t_1,t_2)} \frac{[m(t_1,t_2)]^n}{n!},
	$$
	where $m(t_1,t_2)= 
	{\int_{t_1}^{t_2} \lambda(s)ds}$ and
	$\lambda(t)$ is the arrival rate.
Unlike the Poisson process (where $\lambda(t)=\lambda$), NHPP has nonstationary increments and this makes the use of NHPP suitable for modelling ED arrival process, which is usually strongly time--varying. Of course, appropriate statistical tests must be applied to available data to check if NHPP fits. This is usually performed by assuming that NHPP has a rate which can be considered approximately piecewise constant. Hence, Kolmogorov--Smirnov (KS) statistical test can be applied in separate and equally spaced intervals and usually the classical Conditional--Uniform (CU) property of the Poisson process is exploited  \cite{Brown.2005,Kim.2014,Kim.2014b}. Unlike standard KS test, in the CU~KS test the data are transformed before applying the test. More precisely, by CU property, the piecewise constant NHPP is transformed into a sequence of i.i.d. random variables uniformly distributed on $[0, 1]$ so that it can be considered a (homogeneous) Poisson process in each interval. In this manner, the data from all the intervals can be merged in a single sequence of i.i.d. random variables uniformly distributed on $[0, 1]$. 
This procedure, proposed in \cite{Brown.2005}, enables to remove nuisance parameters and to obtain independence from the rate of the Poisson process on each interval. Hence data from separate intervals (with different rates on each of them) and also from different days can been combined, avoiding common drawback due to large within-day and day-to-day variation of the ED patient arrival rate. Actually, Brown et al. in \cite{Brown.2005} apply CU KS test after performing a further logarithmic data transformation. In \cite{Kim.2014b,Kim.2015}, this approach has been extensively tested along with alternative data transformations proposed in early papers \cite{Durbin.1961} and \cite{Lewis.1965}. However, Kim and Whitt in \cite{Kim.2014} observed that this procedure applied to ED patient arrival data is fair only if they are ``analyzed carefully''. This is due to the fact that the following three issues must be seriously considered: {\em 1) data rounding}, {\em 2) choice of the intervals}, {\em 3) overdispersion}. In fact, the first issue may produce batch arrivals (zero length interarrival times) that are not included in a NHPP, so that unrounded data (or an unrounding procedure) must be considered. The second is a major issue in dealing with ED patient arrivals, since arrival rate can rapidly change so that the piecewise constant approximation is reasonable only if the interval are properly chosen. The third issue regards combining data from multiple days. Indeed, in studying ED patient arrival process, it is common to combine data from the same time slot from different weekdays, being this imperative when data from a single day are not sufficient for statistical testing. Data collected from EDs database usually show large variability over successive weeks mainly due to seasonal phenomena like flu season, holiday season, etc. However, this overdispersion phenomenon must be checked by using a dispersion test on the available data (e.g. \cite{Kathirgamatamby:53}).
\par
In this work we propose a new modelling approach for ED patient arrival process based on a piecewise constant approximation of the arrival rate accomplished with {\em non equally spaced intervals}. This choice is suggested by the typical situation that occurs in EDs where the arrival rate is low and varying during the night hours, and it is higher and more stable in the day time, this is indeed what happens in the chosen case study. Therefore, to obtain an accurate representation of the arrival rate $\lambda(t)$ by a piecewise constant function $\lambda_D(t)$, a finer discretization of the time--domain is required during the night hours, as opposite to day time. For this reason, the proposed method finds the best partition of the 24 hours into intervals not necessarily equally spaced. 
\par
As far as the authors are aware, the use of an optimization method for identifying stochastic processes characterizing the patient flow through an ED was already proposed in \cite{Guo.2016}, but that study aimed at determining the optimal service time distribution parameters (by using a meta--heuristic approach) and it did not involve ED arrival process. Therefore our approach represents the first attempt to adopt an optimization method for determining the best stochastic model for the ED process arrivals. In the previous work \cite{DGLMR:tr} 
a preliminary study was performed following the same approach. Here, with respect to \cite{DGLMR:tr}, we propose a significantly enhanced statistical model which allows us to obtain better results on the case study we consider.
\par
In constructing a statistical model of the ED patient arrivals, a natural way to define a selection criterion is to evaluate the fit error between $\lambda(t)$ and its approximation $\lambda_D(t)$. However, the true arrival rate is unknown. In the approach we propose, as opposite to \cite{Kim.2014}, no analytical model is assumed for $\lambda(t)$, but it is substituted by an ``empirical arrival rate model'' $\lambda_F(t)$ obtained by a sample approximation corresponding to the very fine uniform partition of the $24$ hours into intervals of $15$ minutes. In each of these intervals the average arrival rate values has been estimated from data obtained by collecting samples over the same day of the week, for all the weeks in some months
using experimental data for the ED patient arrival times. 
Hence, any other $\lambda_D(t)$ corresponding to a grosser partition of the day must be compared to $\lambda_F(t)$. In other words, an optimization problem is solved to select the best day partition in non equally spaced intervals, determining a piecewise constant approximation of the arrival rate over the 24 hours with the best fit to the empirical model. Therefore, the objective function (to be minimized) of the optimization problem we formulate, comprises the fit error, namely the mean squared error. Moreover, an additional penalty term is included aiming at obtaining the overall regularity of the optimal approximation, being the latter measured by means of the sum of the squares of the jumps between the values in adjacent intervals. The rationale behind this term is to avoid optimal solutions with too rough behavior, namely few long intervals with high jumps.  
\par  
To make the result reliable, a number of constraints must be considered. First, the length of each interval of the partition can not be less than a fixed value (half an hour, one hour). Moreover, for each interval, 
\begin{itemize}	
\item the CU KS test must be satisfied to support the NHPP hypothesis;
\item the dispersion test must be satisfied to ensure that data are not overdispersed, and could be considered as a realization of the same process (no week seasonal effects).
\end{itemize}
\par
The resulting problem is a black-box constrained optimization problem and to solve it we use a method belonging  to the class of Derivative-Free Optimization. In particular we use the new algorithmic framework recently proposed in \cite{LiuzziLucidiRinaldi:2020} which handles black-box problems with integer variables.  
\par
We performed an extensive experimentation on data collected from the ED of a big hospital in Rome (Italy), also including some significant sensitivity analyses. 
The results obtained show that this approach enables to determine the number of intervals and their length such that an accurate approximation of the empirical arrival rate is achieved, ensuring the consistency between the NHPP hypothesis and the arrival data. The regularity of optimal piecewise constant approximation can be also finely tuned by proper weighing a penalty term in the objective function with respect to the fit error term.
\par
It is worth noting that the use a piecewise constant function for approximating  the arrival rate function is usually required by the most common discrete event simulation software packages  when implementing ED patient arrivals process as a NHPP.
\par	
The paper is organized as follows. In Section~\ref{sec:casestudy}, we briefly report information on the hospital ED under study. Section~\ref{sec:statistical_model} describes the statistical model we propose. The optimization problem we consider is stated in Section~\ref{sec:opt} and the results of an extensive experimentation are reported in Section~\ref{sec:results}. Finally Section~\ref{sec:conclusions} includes some concluding remarks.

\section{The case study under consideration}\label{sec:casestudy}
The case study we consider concerns the ED of the {\em Policlinico Umberto~I}, a very large hospital in Rome, Italy. It is the biggest ED in the Lazio region in terms of yearly patients arrivals (about 140,000 on the average). Thanks to the cooperation of the ED staff, we were able to collect data concerning the patient flow through the ED for the whole year 2018. In particular, for the purpose of this work, we focus on the patient arrivals data collected in the first $m$ weeks of the year. Both walk-in patients and patients transported by emergency medical service vehicles are considered.
\par
\begin{figure}[htbp]
		\centering
		\includegraphics[width=11truecm,height=5truecm]{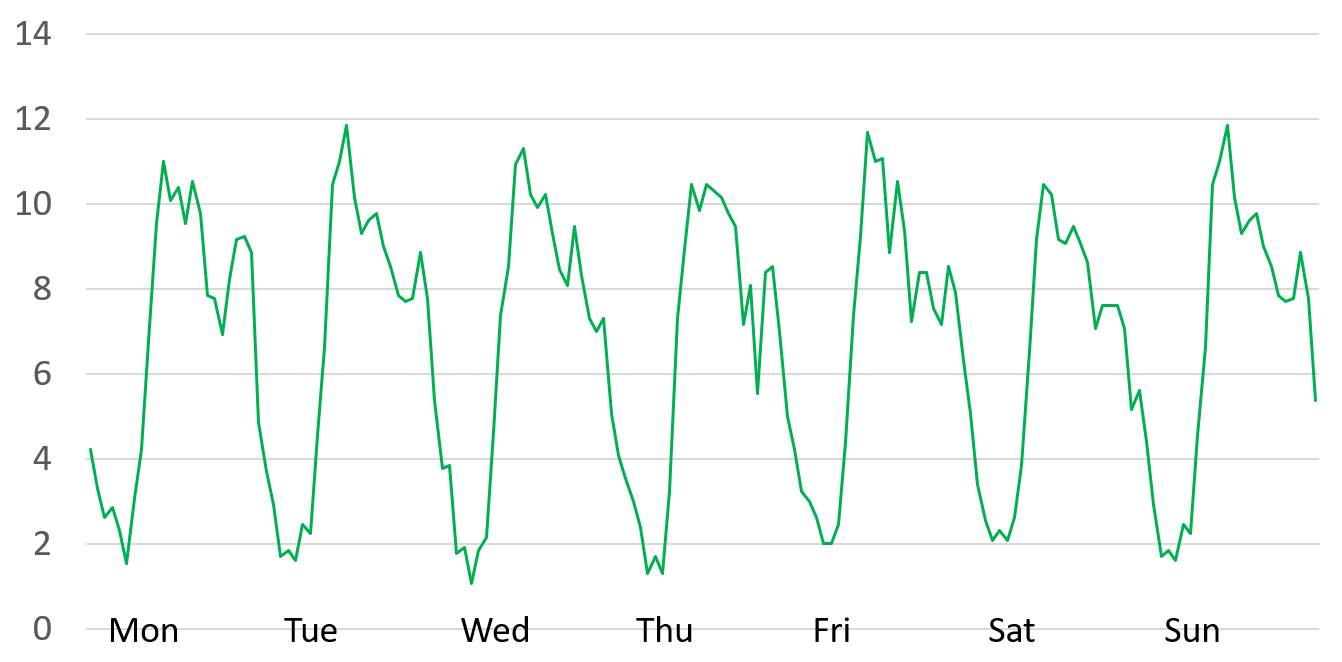}
		\caption{Plot of the weekly average arrival rate for the first 13 weeks of the year.}\label{fig:weekly_arrival_rate}
	\end{figure}

In Figure \ref{fig:weekly_arrival_rate}, the weekly hourly average arrival rate to the ED is shown for $m=13$, i.e. for data collected from the 1st of January to the 31st of March.
In particular, for each day of the week, the arrival rate is obtained by averaging the number of arrivals occurring in the same hourly time slot over the 13 weeks considered. We observe that, in accordance with the literature (see, i.e., \cite{Kim.2014}), the average arrival rates among the days of the week are significantly different. Therefore, since averaging over these days would lead to inaccurate results, the different days of the week must be considered separately.
\par
\begin{figure}[htbp]
	\centering
	\includegraphics[width=11truecm, height=7truecm]{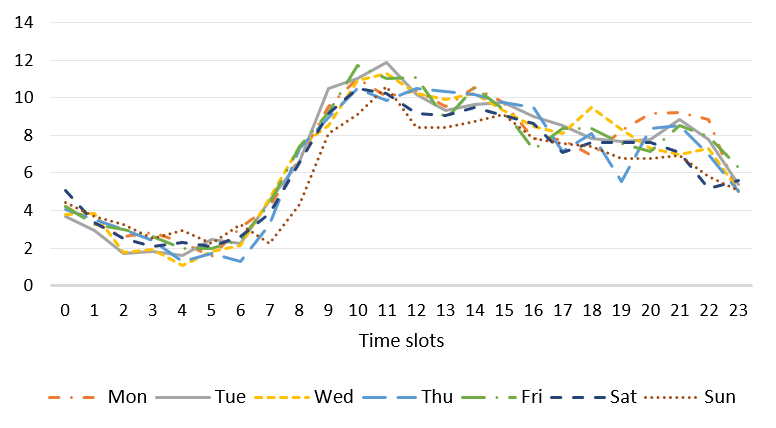}
	\caption{Plot of the comparison among hourly average arrival rate for each day of the week for the first 13 weeks of the year.}
	\label{fig:daily_arrival_rate_comparison}
\end{figure}

 
Figure \ref{fig:daily_arrival_rate_comparison} reports the hourly average arrival rate for each day of the week,
again referring to $m=13$, i.e. to the first 13 weeks of the year.
By observing this figure, we expect that, being the shape of each rate similar, the approach proposed in this work allows to obtain similar partitions of the 24 hours on different days of the week. This enables to focus only on one arbitrary day of the week. Specifically, Tuesday is the day chosen to apply the methodology under study, since the shape of its arrival rate can be considered representative of the other days.      
\par
\begin{figure}[htbp]
	\centering
	\includegraphics[width=11truecm,height=7truecm]{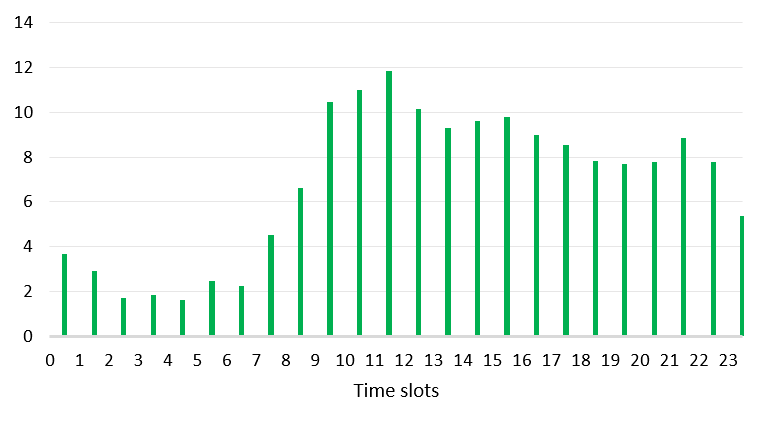}
	\caption{Plot of the average hourly arrival rate for the Tuesdays over the 13 considered weeks of the year.}\label{fig:arrival_rate}
\end{figure}

\begin{figure}
	\centering
	\includegraphics[width=12truecm,height=7truecm]{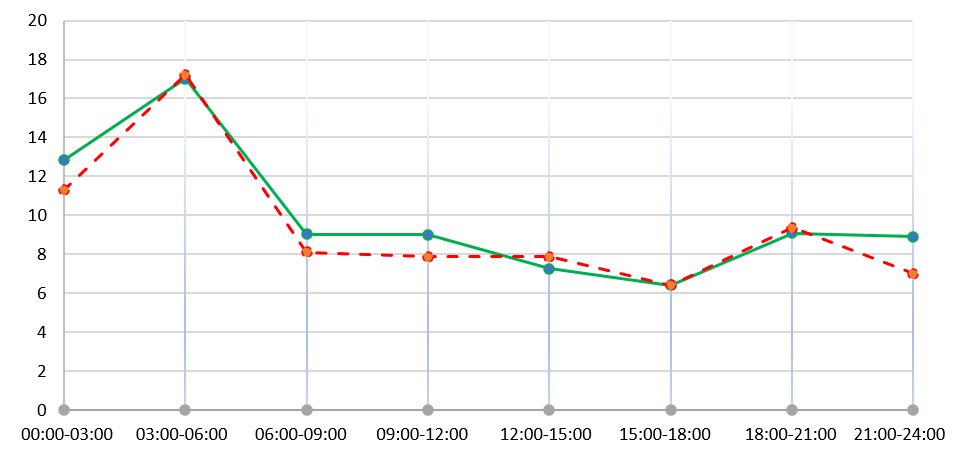}
	\caption{Plot of the average (in solid green) and variance (in dashed red) of the interarrival times
	for the first Tuesday of year 2018. On the abscissa axis, 3-hours time slots are considered.}
	\label{fig:maggio}
\end{figure}

	
In Figure \ref{fig:arrival_rate}, the plot of the hourly average arrival rate for the Tuesdays over the 13 considered weeks is reported, while Figure \ref{fig:maggio} shows mean and variance of the interarrival times occurred on the first Tuesday of the year 2018.
From this latter figure, we observe that these two statistics have similar values within each 3-hours time slot and this is in accordance with the property of the Poisson probability distribution for which mean and variance coincide.

\section{Statistical model}\label{sec:statistical_model}
The arrival process at EDs is usually characterized by a strong within-day variation both in the arrival rate and interarrival times: experimental data show rapid changes in the number of arrivals during the night hours, as opposite to a smoother profile at day time.
As we already mentioned in the Introduction, for this reason the ED arrival process is usually modeled as a NHPP. 
\par
No analytical model is available for the arrival rate $\lambda(t)$, and therefore a suitable representation of the unknown function is needed. A realistic representation can be obtained by averaging the number of arrivals observed in experimental data on suitable intervals over the 24 hours of the day, non necessarily equally spaced. 
Let $\{T_i\}$ denote a partition $P$ of the observation period $T = [0,24]$ (hours) in $N$ intervals, and let $\{\lambda_i\}$ be corresponding sample average rates. Then a piecewise constant approximation of $\lambda(t)$ is written as follows
\begin{equation}
\lambda_D(t)= \sum_{i=1}^N \lambda_i\,\textbf{1}_{T_i}(t), \quad t\in T
\end{equation}
where $\textbf{1}_{T_i}(t)$ is 1 for $t\in T_i$ and 0 otherwise (the indicator function of set $T_i$). 
Any partiton $P$ gives rise to a different approximation $\lambda_D(t)$, depending on the number of intervals and their lengths. Therefore a criterion is needed to select the best partition $P^\star$ with some desirable features. 
\par
First of all, we need to ensure that there is no overdispersion in the arrival data. We refer to the commonly used {\em dispersion test} proposed in \cite{Kathirgamatamby:53} and reported in \cite{Kim.2014}. If it is satisfied, then it is possible to combine arrivals for the same day of the week over different weeks. To this aim, for any partition $P$, let $\{k_i^r\}$ denote the number of arrivals in the $i$-th partition interval $T_i$ in the $r$-th week, $ r=1,\ldots, m$.
Consider the statistics 
$$ Ds_i = \displaystyle\frac{1}{\mu_i}\displaystyle\sum_{r=1}^m \left( k_i^r - \mu_i\right)^2, \quad i=1,\ldots,N,
$$
where $\mu_i = \frac{1}{m}\sum_{r=1}^m k_i^r$
is the average number of arrivals in the given interval for the same day of the week over the considered $m$ weeks. Under the null hypothesis that the counts $\{k_i^r\}$ 
are a sample of $m$  independent Poisson random variables with the same mean count $\mu_i$ (no overdispersion), then $Ds_i$ is distributed as $\chi^2_{m-1}$, the chi-squared distribution with $m-1$ degrees of freedom. Therefore the null hyphotesis is accepted with $1-\alpha$ confidence level if
\begin{equation}\label{eq:disptest}
Ds_i \leq \chi^2_{m-1,\alpha}, \quad i=1,\ldots,N,
\end{equation}
where $\chi^2_{m-1,\alpha}$ is of course the $\alpha$ level critical value of the $\chi^2_{m-1}$ distribution. 
\par
Furthermore, the partition is {\em feasible} if data are consistent with NHPP. Namely, if we denote by $k_i$ the number of arrivals in each interval $T_i=[a_i, b_i)$ obtained by considering data of the same weekday, in the same interval, over $m$ weeks, i.e. $k_i=\sum_{r=1}^m k_i^r$, $i=1,\ldots,N$, the partition is feasible if each $k_i$ has a Poisson
distribution with rate $\lambda_i$ obtained as $\mu_i/(b_i-a_i)$. 
To check the validity of the Poisson hypothesis, the CU KS test can be performed (see \cite{Brown.2005,Kim.2014}). We prefer to use CU KS with respect to Lewis KS test since this latter is highly sensitive to rounding of the data and moreover CU KS test has more power against alternative hypotheses involving exponential interarrival times (see \cite{Kim.2014b} for a detailed comparison between the effectiveness of the two tests).
\par
To perform CU KS test, for any interval $T_i=[a_i, b_i)$, let $t_{ij}$, $j=1, \ldots ,k_i,$ be the arrival times within the $i$-th interval obtained as union over the $m$ weeks of the arrival times in each $T_i$. Now consider the rescaled arrival times defined by $\tau_{ij} =\displaystyle\frac{t_{ij}-a_i}{b_i-a_i}$. The rescaled arrival times, conditionally to the value $k_i$, are a collection of i.i.d. random variables uniformly distributed over $[0,1]$. Hence, in any interval we compare the theoretical cumulative distribution function (cdf) $F(t) = t$ with the empirical cdf 
$$F_i(t) = \frac{1}{k_i}\sum_{j=1}^{k_i} \textbf{1}_{\{\tau_{ij} \le t\}}, \qquad 0 \le t \le 1.$$
The test statistics is defined as follows
\begin{equation} \label{D_i}
D_i = \sup_{0 \le t \le 1}(\vert F_i(t)-t\vert).
\end{equation}
The critical value for this test is denoted as $T(k_i,\alpha)$ and its values can be found on the KS test critical values table. Accordingly, the Poisson hypothesis is accepted if 
\begin{equation}\label{eq:testKS}
D_i \le T(k_i,\alpha), \quad i=1, \ldots , N.
\end{equation}
This test has to be satisfied on each interval $T_i$ to qualify the partition $P$ given by $\{T_i\}$ as feasible, in the sense that CU KS test is satisfied, too.
\par
A further restriction is imposed on the feasible partitions. Given the experimental data, realistic partitions can not have a granularity too fine to avoid that some $k_i$ being too small may unduly determine the rejection of the CU KS test. To this aim the value of 1 hour was chosen as lower threshold value, taking into account the specific case study considered (see also Figure \ref{fig:arrival_rate}).
\par
\begin{figure}[t]
	\centering
	\includegraphics[width=11truecm,height=7truecm]{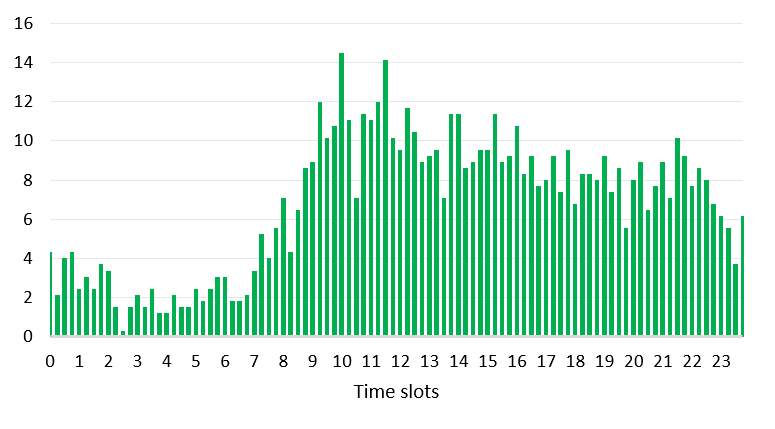}
	\caption{Plot of the daily average arrival rate $\lambda_F(t)$ with intervals of 15 minutes.}\label{fig:daily_arrival_rate15min}
\end{figure}

Now let us evaluate the feasible partitions also in terms of the characteristics of function $\lambda_D(t)$. It would be amenable to define a fit error with respect to $\lambda(t)$, which unfortunately is unknown. The problem can be got around by considering a piecewise constant approximation $\lambda_F(t)$ over a very fine partition ${P}_F$ of $T$. A set of 96 equally space intervals of $15$ minutes was considered and the corresponding average rates $\lambda_i^F$ were estimated from data. The plot of $\lambda_i^F$ is reported in Figure \ref{fig:daily_arrival_rate15min}.
The function $\lambda_F(t)$ can be considered as an {\em empirical arrival rate model}. Note that partition ${P}_F$ need not be feasible since it only serves to define the finest piecewise constant approximation of $\lambda(t)$. Therefore the following fit error can be defined 
\begin{equation} \label{Fit}
E(P)= \sum_{j=1}^N \sum_{i_j=1}^{N_j} (\lambda_j-\lambda_{i_j}^F)^2
\end{equation}
where $N_j$ is the number of intervals of $15$ minutes contained in $T_j$, and identified by the set of indexes $\{i_j\}\subset\{1,\ldots,96\}$. 
\par
Finally it is also advisable to characterize the ``smoothness'' of any approximation $\lambda_D(t)$ to avoid very gross partitions with high jumps between adjacent intervals
by means of the mean squared error
\begin{equation}\label{Smooth}
S(P) = \sum_{j=2}^N (\lambda_j-\lambda_{j-1})^2.
\end{equation}
\par
In the following Section~\ref{sec:opt} the model features illustrated above are organized in a proper optimization procedure that provides the selection of the best partition according to conflicting goals. 
\par
The approach we propose enables to well address the major two issues raised in \cite{Kim.2014} (and reported in 
the Introduction) when dealing with modelling ED patient arrivals, namely the {\em  choice of the intervals} and the {\em  overdispersion}. As concerns the third issue, the {\em data rounding}, the arrival times in the data we collected are rounded to seconds (format {\tt hh:mm:ss}), and actually occurrences of simultaneous arrivals which would cause zero interarrival times are not present. Therefore, we do not need any unrounding procedure. Anyhow, as already pointed out above, the CU~KS test we use is not very sensitive to data rounding.

\section{Statement of the optimization problem}\label{sec:opt}
Any partition $P=\{T_i\}$ of $T = [0,24]$ is characterized by the boundary points $\{x_i\}$ of its intervals and by their number $N$. Let us introduce a vector of variables $x\in \mathbb{Z}^{25}$ such that
$$T_i=[x_i, x_{i+1}),$$
$i=1, \ldots ,24$, with $x_1=0$ and $x_{25}=24$.
\par
Functions in (\ref{Fit}) and (\ref{Smooth}) are indeed functions of $x$, and therefore will be denoted by $E(x)$ and $S(x)$, respectively. Therefore, the objective function that constitutes the selection criterion is given by
\begin{equation}\label{ObjPen}
f(x) = E(x) + w S(x),
\end{equation}
where $w>0$ is a parameter that controls the weight of the smoothness penalty term with respect to the fit error: the larger $w$, the smaller the difference between average arrival rates in adjacent intervals; this in turn implies that on a steep section of $\lambda_F(t)$ an increased number of shorter intervals is adopted to fill the gap with relatively small jumps.

The set $\cal P$ of feasible partitions is defined as follows:
\begin{equation} \label{Feasibleset}
\begin{array}{l}
{\cal P}=\Bigl\{x\in{\mathbb{Z}}^{25} ~ | ~ x_1 = 0, \quad x_{25}=24, \quad x_{i+1}-x_i\geq \ell_i, \quad g_i(x)\leq 0, \bigr. \cr
\ \cr
\bigl.\hspace{1truecm}   h_i(x)\leq 0,  \quad i=1, \ldots , N \Bigr\}
\end{array}
\end{equation}
where
\begin{eqnarray} 
\ell_i &=&
\begin{cases}
0 \quad \hbox{if} \quad x_i=x_{i+1},\\
\ell \quad \hbox{otherwise},\\
\end{cases} \label{eq:vincolibound} \\
\ \nonumber \\
g_i(x) &=&
\begin{cases}
0 \quad \hbox{if} \quad x_i=x_{i+1},\\
D_i - T(k_i,\alpha) \quad \hbox{otherwise},\\
\end{cases} \label{eq:vincolistrong}  \\
\ \nonumber \\
h_i(x) &=&
\begin{cases}
0 \quad \hbox{if} \quad x_i=x_{i+1},\\
Ds_i - \chi^2_{m-1,\alpha}  \quad \hbox{otherwise},\\
\end{cases} \label{eq:vincoliweak} \\ \nonumber
\end{eqnarray}
$i=1, \ldots , N$. The value $\ell$  in
(\ref{eq:vincolibound})	denotes the minimum interval length allowed and we assume $\ell\geq 1/4$. Of course, constraints $g_i(x)\leq 0$ represents the satisfaction of the CU~KS test in (\ref{eq:testKS}), while constraints $h_i(x)\leq 0$ concern the dispersion test in (\ref{eq:disptest}). 
Therefore, the best piecewise constant approximation $\lambda_D^\star(t)$ of the time-varying arrival rate $\lambda(t)$ is obtained by solving the following black-box optimization problem:
\begin{equation}\label{optBestFit}
\begin{split}
\max ~~ & f(x) \\
s.t. ~~ & x\in {\cal P}. \\
\end{split}
\end{equation}
We highlight that the idea to use as constraints of the optimization problem a test to validate the underlying statistical hypothesis on data along with a dispersion test
is completely novel in the framework of modeling ED patient arrivals process. The only proposal which use a similar approach is in our previous paper \cite{DGLMR:tr}.
\par
It is important to note that in (\ref{ObjPen}) the objective function has not analytical structure with respect to the independent variables and it can only be computed by a data-driven procedure once the $x_i$'s values are given. The same is true for the constraints $g_i(x)$ and $h_i(x)$ in (\ref{Feasibleset}), too. Therefore the problem in hand is an integer nonlinear constrained black-box problem, and both the objective function and the constraints are relatively expensive to compute and this makes it difficult to efficiently solve. In fact, classical optimization methods either can not be applied (since based on the analytic knowledge of the functions involved) or they are not efficient especially when evaluating the functions at a given point is very computationally expensive. Therefore to tackle problem (\ref{optBestFit}) we turned our attention to the class of Derivative-Free Optimization and black-box methods (see, e.g.,  \cite{Audet.2017,Conn.2009,Larson.2019}). More in particular, we adopt the algorithmic framework recently proposed in \cite{LiuzziLucidiRinaldi:2020}. It represents a novel strategy for solving black-box problems with integer variables and it is based on the use of suited search directions and a nonmonotone linesearch procedure. Moreover, it can handle generally-constrained problems by using a penalty approach. 
We refer to \cite{LiuzziLucidiRinaldi:2020} for a detailed description and we only highlight that the results reported in \cite{LiuzziLucidiRinaldi:2020} clearly show that this algorithm framework is particularly efficient in tackling black-box problems like the one in (\ref{optBestFit}). In particular, the effectiveness of the adopted exploration strategy with respect to state-of-the-art methods for black-box is shown. This is due to the fact that the approach proposed in \cite{LiuzziLucidiRinaldi:2020} combines computational efficiency with a high level of reliability.

\section{Experimental results}\label{sec:results}
In this section we report the results of an extensive experimentation on data concerning the case study described in Section~\ref{sec:casestudy}, namely the ED patient arrivals collected in the first $m$ weeks of year 2018. Different values of the number of weeks $m$ have been considered. Standard significance level $\alpha=0.05$ is used the CU KS and dispersion tests. 
\par
As regards the optimization problem in hand the value of $\ell$ in (\ref{eq:vincolibound}) is set to $1$ hour. Moreover, it is important to note that different values of the weight $w$ in the objective function (\ref{ObjPen}) lead to various piecewise constant approximations with a different fitting accuracy and degree of regularity. Therefore, we performed a careful tuning of this parameter, aiming at determining a value which represents a good trade-off between a small fit error and the smoothness of the approximation.
\par
As concerns the parameter values of the optimization algorithm used in our experimentation, we used the default ones (see \cite{LiuzziLucidiRinaldi:2020}). The stopping criterion is based on the maximum number of function evaluations set to 5000. As starting point $x^0$ of the optimization algorithm we adopt the following  
\begin{equation}\label{eq:startingpoint}
x^0_i=i-1, \qquad i=1,\ldots,25,
\end{equation}
which corresponds to the case of 24 intervals of unitary length. This choice is a commonly used partition in most of the approaches proposed in literature (see e.g. \cite{Ahalt.2018,Kim.2014}).
Table~\ref{tab:CUKSTest-initial} in the Appendix reports the results of CU KS and dispersion tests applied to the partition corresponding to the starting point $x^0$, considering $m=13$ weeks. In particular, in Table~\ref{tab:CUKSTest-initial} for each one-hour slot the sample size $k_i$ is reported along with the $p$-value and the acceptance/rejection of the  null hypothesis of the corresponding test.
We observe that the arrivals are not overdispersed in any interval of the partition corresponding to $x^0$, i.e. all the constraints $h_i(x)\leq 0$ are satisfied and this allows us to combine data for the same day of the week over successive weeks. However, this  partition is even unfeasible, i.e. $g_i(x) >0$ for some $i$; this corresponds to reject the statistical hypothesis on some $T_i$. Notwithstanding, even if the starting point is unfeasible, the optimization algorithm we use is able find a feasible solution which minimizes the objective function.
\par
As we already mentioned, the choice of a proper value for the weight $w$ in the objective function (\ref{ObjPen}) is important and not straightforward. On the other hand, the number $m$ of the considered weeks also affects both the accuracy of the approximation, through the average rates estimated on each interval, and the consistency of the results, which is ensured by constraints \eqref{eq:vincolistrong} and \eqref{eq:vincoliweak}. However, while $w$ is related to the statement of the optimization problem \eqref{optBestFit} and it can be arbitrarily chosen, the choice of $m$ is strictly connected to the available data. In \cite[Section~4]{Kim.2014}, the authors  assert that, having 10 arrivals in the one-hour slot 9--10 a.m., it is necessary to combine data over 20 weeks in order to have a sufficient sample size (200 patient arrivals). However, being their approach based on equally-spaced intervals, one-hour slots are also adopted during off-peak hours, for instance during the night. This implies that the sample size corresponding to data combination over 20 weeks for these  slots could no longer be sufficient to guarantee good results. This is clearly pointed out in Table~\ref{tab:CUKSTest-initial} where the sample size $k_i$ corresponding to some of the one-hour night slots is very low considering $m=13$ weeks and it remains insufficient even if 26 weeks are considered (see subsequent Table~\ref{tab:Test-initial}). The approach we propose overcomes this drawback since, for each choice of $m$,  we determine the length of the intervals as solution of the optimization problem \eqref{optBestFit}. Of course, there could be values of $m$ such that problem \eqref{optBestFit} has not feasible solutions, i.e. a partition such that the NHPP hypothesis holds and the results are consistent does not exists for such $m$.
\par
In order to deeper examine how the parameters $w$ and $m$ affect the optimal partition, we performed a sensitivity analysis, focusing first on the case with fixed $m$ and $w$ varying.
In particular, we have chosen to focus on $m = 13$ weeks, which enables to achieve an optimal solution by running the optimization algorithm without overly computational burden. Anyhow, we expect that no substantial changes in the conclusions would be obtained with different values of $m$ and this is confirmed by further experimentation whose results are not reported here for the sake of brevity. 
\par
This analysis allows us to obtain several partitions that may be considered for a proper fine-tuning of $w$. In particular, we consider different values of $w$ within the set $\{0, ~ 0.1,~ 1, ~ 10, ~ 10^3\}$. Table~\ref{tab:CUKSTest-optimal} in the Appendix reports the optimal partitions obtained by solving problem \eqref{optBestFit} for these values of $w$. In particular, Table~\ref{tab:CUKSTest-optimal} includes the intervals of the partition, the value of the sample size $k_i$ corresponding to each interval over $13$ weeks and the results of the CU KS and dispersion tests, namely the $p$-value and the acceptance/rejection of the null hypothesis of the corresponding test.
\par
\begin{figure}[htbp]
	\caption{Graphical comparison between the empirical arrival rate model $\lambda_F(t)$ (in green) and the piecewise constant approximation $\lambda_D(t)$ (in red) corresponding to the optimal partition obtained by solving problem \eqref{optBestFit} for different values of the parameter $w$. From top to bottom: $w=0, 0.1, 1, 10, 10^3$.}\label{fig:grafici-partizioni}
	\centering
	\includegraphics[width=8truecm]{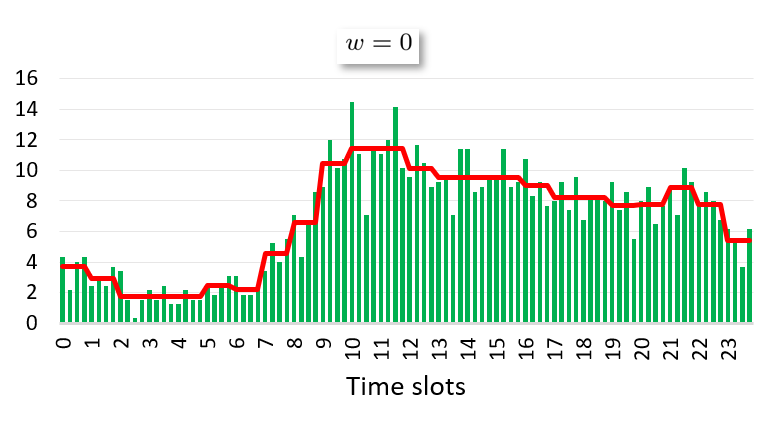}
	\quad
	\includegraphics[width=8truecm]{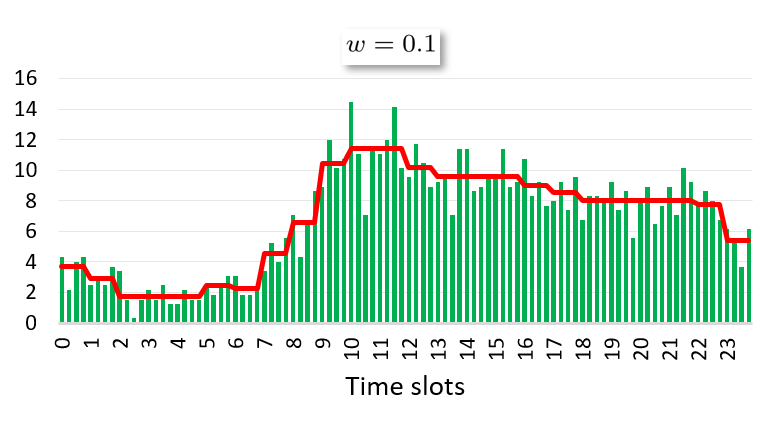}
	\quad
	\includegraphics[width=8truecm]{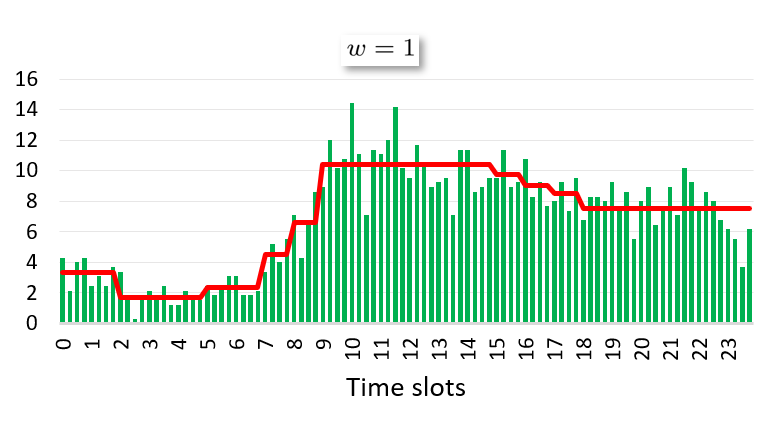}
	\quad
	\includegraphics[width=8truecm]{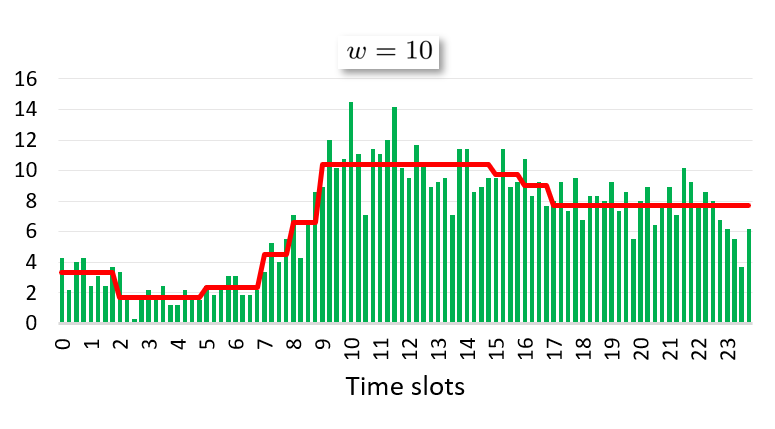}
	\quad
	\includegraphics[width=8truecm]{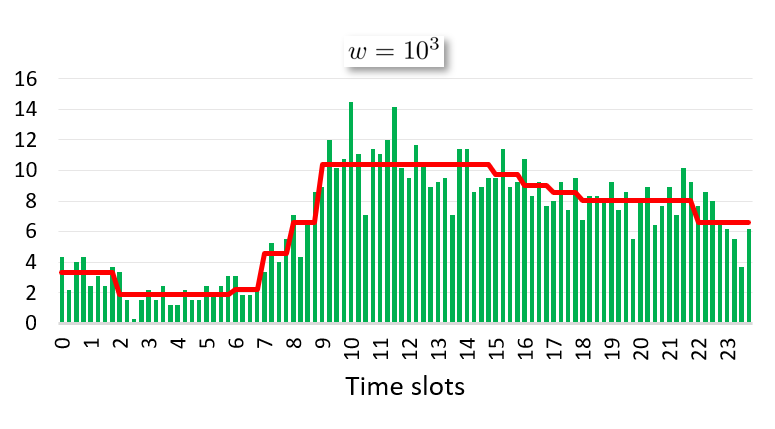}
\end{figure}

In Figure~\ref{fig:grafici-partizioni}, for graphical comparison,  we report the plots of the empirical arrival rate model $\lambda_F(t)$ and its piecewise constant approximation  $\lambda_D(t)$ corresponding to the optimal partitions obtained.
Two effects can be clearly observed as $w$ increases: on the one hand, on steep sections of $\lambda_F(t)$, shorter intervals are adopted to reduce large gaps between adjacent intervals; on the other hand, when $\lambda_F(t)$ is approximately flat, a lower number of intervals may be sufficient to guarantee small gaps. This is confirmed by the two top plots in Figure~\ref{fig:grafici-partizioni} which correspond to $w=0$ and $w=0.1$. In fact, in the first plot ($w=0$) where only the fit error is included in the objective function and in the second one ($w=0.1$) where anyhow the fit error is the dominant term of the objective function, the optimal partition is composed by a relatively large number of intervals. In particular, in the partition corresponding to $w = 0.1$, fewer intervals are adopted during the day time. As expected, a smaller number of intervals is attained when $w = 1$, $w = 10$ and $w = 10^3$. Note that, since on the steep section corresponding to the time slot 7:00--10:00 a.m. the maximum number of allowed intervals (due to the lower threshold value of one hour given by the choice $\ell=1$ in \eqref{eq:vincolibound}) is already used, the only way to decrease the smoothness term of the objective function is to enlarge the intervals during both the day and the night. It is worth noting that for $w=10^3$, the number of intervals increases if compared with the case $w=10$. This occurs to offset the increase in the fit error term due to the use of a smaller number of intervals on the flatter sections. As a consequence, the partition has an unexpected interval at the end of the day.   
\par
We point out that for each value of $w$, the optimization algorithm finds an optimal partition (of course feasible with respect to all the constraints), despite some constraints related to the CU KS test are violated in the initial partition, i.e. the one corresponding to $x^0$ in \eqref{eq:startingpoint}, namely the standard assumption of one-hour slots usually adopted. This means that the used data are in accordance with the NHPP hypothesis and they are sufficient to appropriately define the piecewise constant approximation of the ED arrival rate.
\par
Conversely, when the optimization algorithm does not find a feasible partition, the CU KS test or the dispersion test related to some $T_i$ are never satisfied. This implies that the process is not conforming to the NHPP hypothesis or that the data are overdispersed. 
This is clearly highlighted by our subsequent experimentation where we set $w=1$, letting $m$ varying within the set $\{5,9,17,22,26\}$.
\par
First, in Table~\ref{tab:Test-initial} in the Appendix we report the results of CU KS and dispersion tests applied to the partition corresponding to the starting point in  $x^0$\eqref{eq:startingpoint},
for these different values of $m$. Once more, this table evidences that the use of equally-spaced intervals of one-hour length during the whole day can be inappropriate. As an example, see the results of the tests on the time slot 02:00--03:00.
Moreover, note that, for all these values of $m$, the initial partition corresponding to the starting point $x_0$ is infeasible, except when $m = 5$. Indeed, the constraints corresponding to CU KS and dispersion tests are violated for some $T_i$, meaning that the validity of the standard assumption of one-hour time slots strongly depends on the time period considered for using the collected data. 
To this aim, a strength of our approach is its ability to assist in the selection of a reasonable value for $m$. 
If there is no value of $m$ such that the optimization algorithm determines an optimal solution (due to unfeasibility), then it may  be inappropriate to consider the ED arrival process in hand as NHPP.
\par
The subsequent Table~\ref{tab:Test-infeasible-final} includes the optimal partitions obtained by solving problem \eqref{optBestFit} for the considered values of $m\in\{5,9,17,22,26\}$. Like the previous tables, Table~\ref{tab:Test-infeasible-final} includes the intervals of the partition, the value of the sample size $k_i$ corresponding to each interval and the results of CU KS and dispersion tests.
For all the considered values of $m$ the optimization algorithm determines an optimal solution with the only exception of $m=26$. In this latter case, the maximum number of function evaluations allowed to the optimization
algorithm is not enough to compute an optimal solution: in fact, 
we obtain an unfeasible solution since the CU KS test related to the last interval of the day is not satisfied. This could be partially unexpected, since more accurate results should be obtained when considering a greater sample size. However, by adding the last four weeks (passing from $m=22$ to $m=26$) which corresponds to the month of June, the data become affected by a seasonal trend and the NHPP assumption is no longer valid.
\par
\begin{figure}[htbp]
	\caption{Graphical comparison between the empirical arrival rate model $\lambda_F(t)$ (in green) and the piecewise constant approximation $\lambda_D(t)$ (in red) corresponding to the optimal partition obtained by solving problem \eqref{optBestFit} for different values of the parameter $m$. From top to bottom: $m=5, 9, 17, 22, 26$.}\label{fig:grafici-partizioni-mesi}
	\centering
	\includegraphics[width=8truecm]{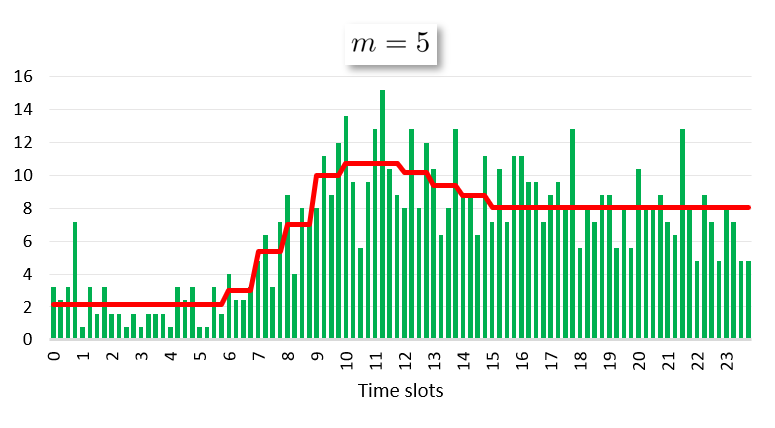}
	\quad
	\includegraphics[width=8truecm]{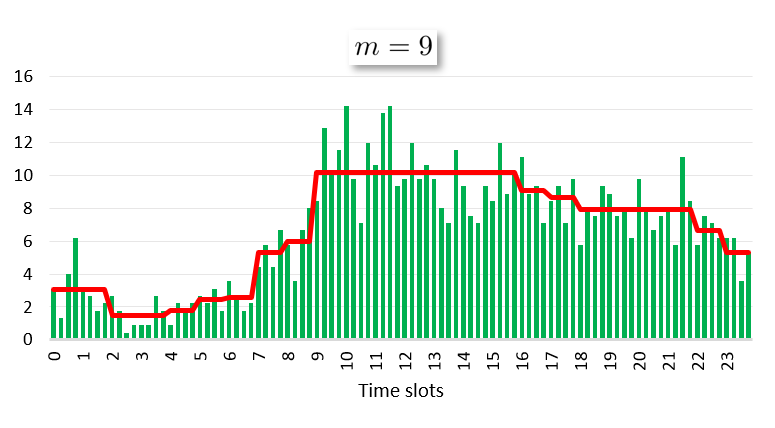}
	\quad
	\includegraphics[width=8truecm]{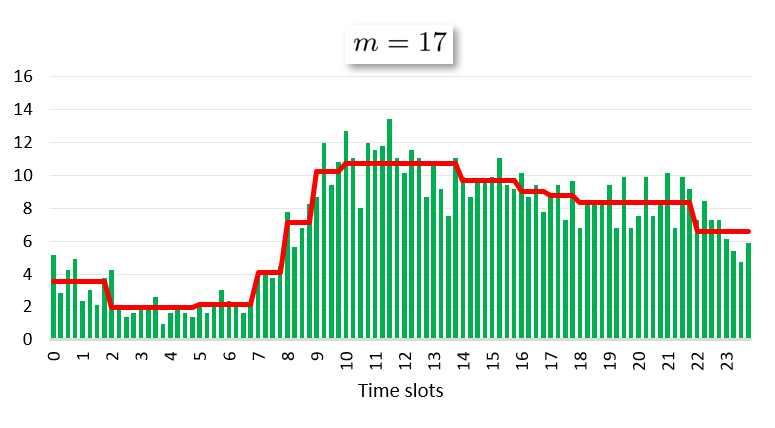}
	\quad
	\includegraphics[width=8truecm]{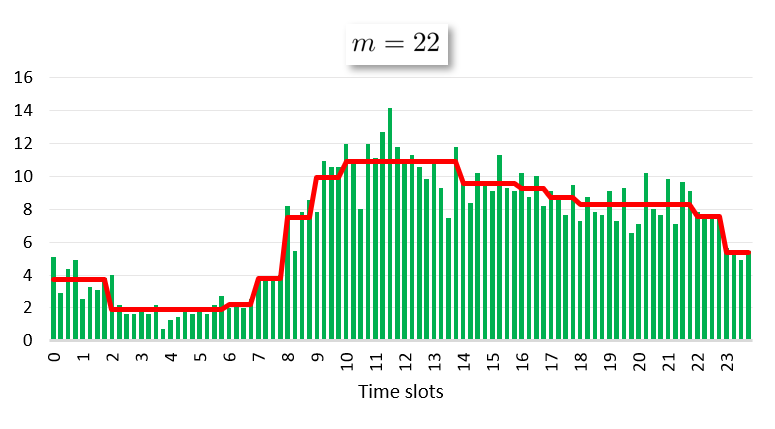}
\end{figure}

In Figure~\ref{fig:grafici-partizioni-mesi} we report a graphical comparison between the empirical arrival rate model $\lambda_F(t)$ and the piecewise constant approximation $\lambda_D(t)$ corresponding to the optimal partitions obtained for the considered values of $m$. We observe that the variability of $\lambda_F(t)$ reduces as the value of $m$ increases since averaging on more data leads to flattening the fluctuation. Despite these rapid oscillations and unlike the other considered values of $m$, for $m = 5$ the empirical model $\lambda_F(t)$ shows a constant trend during both the night and day hours. This results in a piecewise constant approximation $\lambda_D(t)$ that is flat in all the time slots of the 24 hours of the day except the ones related to the morning hours, for which many intervals are used. In fact, to guarantee a good fitting error between $\lambda_D(t)$ and $\lambda_F(t)$, it would be necessary to use shorter intervals, but this is not allowed by the choice $\ell = 1$ in the constraints~\eqref{eq:vincolibound}. For the other considered values of $m$, the number of intervals increases, leading to partitions that improve the fitting error if compared with the case $m = 5$. In particular, we observe that the piecewise constant approximation $\lambda_D(t)$ obtained for $m = 22$ benefits from the lower fluctuations resulting from averaging more data. Therefore, as expected, using the maximum number of available data leads to the most accurate piecewise constant approximation. However, when considering too many data, seasonal phenomena could give rise to the rejection of the null hypothesis of the considered tests, as observed for the case $m = 26$. Moreover, as highlighted at the end of Section 5 in \cite{Kim.2014}, a tendency to reject the NHPP hypothesis (i.e. the null hypothesis of the CU KS test) may be encountered when the sample size is large. In fact, a larger sample size requires a stronger evidence of the null hypothesis in order for the test to be passed. Notwithstanding, our approach is able to overcome these drawbacks, providing us with an optimal strategy to identify the best way of using the collected data.

\section{Conclusions}\label{sec:conclusions}
In this work, we examined the arrival process to EDs by providing a novel methodology that is able to improve the reliability of the modelling approaches frequently used to deal with this complex system, i.e. the Discrete Event Simulation modelling. In accordance with the literature, we adopted the standard assumption of representing the ED arrival process as a NHPP, which is suitable for modelling strongly time-varying processes. In particular, the final goal of the proposed approach is to accurately estimate the unknown arrival rate, i.e. the time-dependent parameter of the NHPP, by using a reasonable piecewise constant approximation. To this aim, an integer nonlinear black-box optimization problem is solved to determine the optimal partition of the 24 hours into a suitable number of non equally spaced intervals. To guarantee the reliability of this estimation procedure, two types of statistical tests are considered as constraints for each interval of any candidate partition: the CU KS test must be satisfied to ensure the consistency between the NHPP hypothesis and the ED arrivals; the dispersion test must be satisfied to avoid the overdispersion of data. To the best of our knowledge, our methodology represents the first optimization-based approach adopted for determining the best stochastic model for the ED arrival process.
\par
The extensive experimentation we performed on data collected from an ED of a big hospital in Italy, shows that our approach is able to find a piecewise constant approximation which represents a good trade-off between a small fit error with the empirical arrival rate model and the smoothness of the approximation. This result is accomplished by the optimization algorithm, despite some constraints in the starting point, which corresponds to the commonly adopted partition composed by one-hour time slots, are violated. Moreover, some significant sensitivity analyses are performed to investigate the fine-tuning of the two parameters affecting the quality of the piecewise constant approximation: the weight of the smoothness of the approximation in the objective function (with respect to the fit error) and the number of weeks considered from the arrivals data. While the former can be arbitrarily chosen by a user according to the desired level of smoothness, the latter affects the accuracy of the arrival rate estimation. In general, the more weeks are considered, the more accurate is the arrival rate approximation, as long as the NHPP assumption still holds and the data do not become overdispersed.
\par
As regards future work, in order to deeper analyze the robustness of the proposed approach, we could use alternative statistical tests, such as the Lewis and the Log tests described in \cite{Kim.2014}, in place of the CU KS test. Moreover, whenever Discrete Event Simulation modelling is the chosen methodology to study ED operation, a model calibration approach could be also used to determine the best value of the weight used in the objective function to penalize the ``smoothness term''. In fact, the optimal value of this parameter could be obtained by minimizing the deviation between the simulation outputs and the corresponding key performance indicators computed through the data. This enables to obtain a representation of the ED arrival process that leads to an improved simulation model of the system under study.

\appendix
\section{Appendix}
In this Appendix we report the detailed results of the CU KS and dispersion tests related to the partitions considered throughout the paper. 
\begin{longtable}{|c|c|c|c|c|c|}
	\caption{Results of the CU KS and dispersion tests (with a significance level of 0.05) applied to each interval of the  partition corresponding to the starting point $x^0$. The considered number of weeks is $m = 13$. For each interval of each partition, the sample size of the dispersion test is $m$. $H_0$ denotes the null hypothesis of the corresponding test.}\label{tab:CUKSTest-initial}\\
	\cline{3-6}
	\multicolumn{2}{c}{} & \multicolumn{2}{|c|}{CU KS test} & \multicolumn{2}{|c|}{Dispersion test}\\ \hline
	Interval & $k_i$ & $p$-value & $H_0$ & $p$-value & $H_0$\\ \hline
	00:00 -- 01:00 & $48$ & $0.836$ & accepted & $0.801$ & accepted \\
	01:00 -- 02:00 & $38$ & $0.950$ & accepted & $0.450$ & accepted \\
	02:00 -- 03:00 & $22$ & $ 0.027$ & rejected & $0.521$ & accepted \\
	03:00 -- 04:00 & $24$ & $0.752$ & accepted & $0.652$ & accepted \\
	04:00 -- 05:00 & $21$ & $0.668$ & accepted & $0.366$ & accepted\\
	05:00 -- 06:00 & $32$ & $0.312$ & accepted & $0.524$ & accepted\\
	06:00 -- 07:00 & $29$ & $0.634$ & accepted & $0.538$ & accepted\\
	07:00 -- 08:00 & $59$ & $0.424$ & accepted & $0.252$ & accepted\\
	08:00 -- 09:00 & $86$ & $0.393$ & accepted & $0.734$ & accepted\\
	09:00 -- 10:00 & $136$ & $0.635$ & accepted & $0.803$ & accepted\\
	10:00 -- 11:00 & $143$ & $0.039$ & rejected & $0.966$ & accepted\\ 
	11:00 -- 12:00 & $154$ & $0.325$ & accepted & $0.999$ & accepted\\ 
	12:00 -- 13:00 & $132$ & $0.858$ & accepted & $0.948$ & accepted\\ 
	13:00 -- 14:00 & $121$ & $0.738$ & accepted & $0.984$ & accepted\\ 
	14:00 -- 15:00 & $125$ & $0.885$ & accepted & $0.500$ & accepted\\ 
	15:00 -- 16:00 & $127$ & $0.928$ & accepted & $0.610$ & accepted\\ 
	16:00 -- 17:00 & $117$ & $0.479$ & accepted & $0.987$ & accepted\\ 
	17:00 -- 18:00 & $111$ & $0.769$ & accepted & $0.516$ & accepted\\ 
	18:00 -- 19:00 & $102$ & $0.458$ & accepted & $0.912$ & accepted\\ 
	19:00 -- 20:00 & $100$ & $0.095$ & accepted & $0.527$ & accepted\\ 
	20:00 -- 21:00 & $101$ & $0.656$ & accepted & $0.586$ & accepted\\ 
	21:00 -- 22:00 & $115$ & $0.763$ & accepted & $0.604$ & accepted\\ 
	22:00 -- 23:00 & $101$ & $0.916$ & accepted & $0.305$ & accepted\\ 
	23:00 -- 24:00 & $70$ & $0.864$ & accepted & $0.104$ & accepted\\ \hline
\end{longtable}
\newpage
\begin{longtable}{|c|c|c|c|c|c|c|}
	\caption{Results of the CU KS and dispersion tests (with a significance level of 0.05) applied to each interval of the optimal partition obtained by solving problem \eqref{optBestFit} for different values of the parameter $w$, with $m$ fixed to 13 weeks. From top to bottom: $w=0,0.1,1,10,10^3$. For each interval of each partition, the sample size of the dispersion test is equal to $m$. $H_0$ denotes the null hypothesis of the corresponding test.}\label{tab:CUKSTest-optimal}\\
	\cline{4-7}
	\multicolumn{3}{c}{} & \multicolumn{2}{|c|}{CU KS test} & \multicolumn{2}{c|}{Dispersion test} \\
	\cline{1-7}
	$w$ & Interval & $k_i$ & $p$-value & $H_0$ & $p$-value & $H_0$ \\ \hline
	& 00:00 -- 01:00 & $48$ & $0.836$ & accepted & $0.801$ & accepted \\
	& 01:00 -- 02:00 & $38$ & $0.950$ & accepted & $0.450$ & accepted \\
	& 02:00 -- 05:00 & $67$ & $ 0.504$ & accepted & $0.100$ & accepted \\
	& 05:00 -- 06:00 & $32$ & $0.312$ & accepted & $0.524$ & accepted\\
	& 06:00 -- 07:00 & $29$ & $0.634$ & accepted & $0.538$ & accepted\\
	& 07:00 -- 08:00 & $59$ & $0.424$ & accepted & $0.252$ & accepted \\
	& 08:00 -- 09:00 & $86$ & $0.393$ & accepted & $0.734$ & accepted\\
	& 09:00 -- 10:00 & $136$ & $0.635$ & accepted & $0.803$ & accepted\\
	0 & 10:00 -- 12:00 & $297$ & $0.433$ & accepted & $0.994$ & accepted\\
	& 12:00 -- 13:00 & $132$ & $0.858$ & accepted & $0.948$ & accepted\\
	& 13:00 -- 16:00 & $373$ & $0.958$ & accepted & $0.502$ & accepted\\
	& 16:00 -- 17:00 & $117$ & $0.479$ & accepted & $0.987$ & accepted\\
	& 17:00 -- 19:00 & $213$ & $0.999$ & accepted & $0.937$ & accepted\\
	& 19:00 -- 20:00 & $100$ & $0.095$ & accepted & $0.527$ & accepted\\
	& 20:00 -- 21:00 & $101$ & $0.656$ & accepted & $0.586$ & accepted\\
	& 21:00 -- 22:00 & $115$ & $0.763$ & accepted & $0.604$ & accepted\\
	& 22:00 -- 23:00 & $101$ & $0.916$ & accepted & $0.305$ & accepted\\
	& 23:00 -- 24:00 & $70$ & $0.864$ & accepted & $0.104$ & accepted\\
	\hline
	& 00:00 -- 01:00 & $48$ & $0.836$ & accepted & $0.801$ & accepted\\
	& 01:00 -- 02:00 & $38$ & $0.950$ & accepted  & $0.450$ & accepted\\
	& 02:00 -- 05:00 & $67$ & $ 0.504$ & accepted & $0.100$ & accepted\\
	& 05:00 -- 06:00 & $32$ & $0.312$ & accepted & $0.524$ & accepted\\
	& 06:00 -- 07:00 & $29$ & $0.634$ & accepted & $0.538$ & accepted\\
	& 07:00 -- 08:00 & $59$ & $0.424$ & accepted & $0.252$ & accepted\\
	& 08:00 -- 09:00 & $86$ & $0.393$ & accepted & $0.734$ & accepted\\
	& 09:00 -- 10:00 & $136$ & $0.635$ & accepted & $0.803$ & accepted\\
	0.1& 10:00 -- 12:00 & $297$ & $0.433$ & accepted & $0.994$ & accepted\\
	& 12:00 -- 13:00 & $132$ & $0.858$ & accepted & $0.948$ & accepted\\
	& 13:00 -- 16:00 & $373$ & $0.958$ & accepted & $0.502$ & accepted\\
	& 16:00 -- 17:00 & $117$ & $0.479$ & accepted & $0.987$ & accepted\\
	& 17:00 -- 18:00 & $111$ & $0.769$ & accepted & $0.516$ & accepted\\
	& 18:00 -- 22:00 & $418$ & $0.660$ & accepted & $0.987$ & accepted\\
	& 22:00 -- 23:00 & $101$ & $0.916$ & accepted & $0.305$ & accepted\\
	& 23:00 -- 24:00 & $70$ & $0.864$ & accepted & $0.104$ & accepted\\
	\hline
	& 00:00 -- 02:00 & $86$ & $0.825$ & accepted & $0.709$ & accepted\\
	& 02:00 -- 05:00 & $67$ & $0.504$ & accepted & $0.100$ & accepted\\
	& 05:00 -- 07:00 & $61$ & $ 0.739$ & accepted & $0.313$ & accepted\\
	& 07:00 -- 08:00 & $59$ & $0.424$ & accepted & $0.252$ & accepted\\
	1 & 08:00 -- 09:00 & $86$ & $0.393$ & accepted & $0.734$ & accepted\\
	& 09:00 -- 15:00 & $811$ & $0.073$ & accepted & $0.955$ & accepted\\
	& 15:00 -- 16:00 & $127$ & $0.928$ & accepted & $0.610$ & accepted\\
	& 16:00 -- 17:00 & $117$ & $0.479$ & accepted & $0.987$ & accepted\\
	& 17:00 -- 18:00 & $111$ & $0.769$ & accepted & $0.516$ & accepted\\
	& 18:00 -- 24:00 & $589$ & $0.059$ & accepted & $0.922$ & accepted\\
	\hline
	& 00:00 -- 02:00 & $86$ & $0.825$ & accepted & $0.709$ & accepted\\
	& 02:00 -- 05:00 & $67$ & $0.504$ & accepted & $0.100$ & accepted\\
	& 05:00 -- 07:00 & $61$ & $ 0.739$ & accepted & $0.313$ & accepted\\
	& 07:00 -- 08:00 & $59$ & $0.424$ & accepted & $0.252$ & accepted\\
	10 & 08:00 -- 09:00 & $86$ & $0.393$ & accepted & $0.734$ & accepted\\
	& 09:00 -- 15:00 & $811$ & $0.073$ & accepted & $0.955$ & accepted\\
	& 15:00 -- 16:00 & $127$ & $0.928$ & accepted & $0.610$ & accepted\\
	& 16:00 -- 17:00 & $117$ & $0.479$ & accepted & $0.987$ & accepted\\
	& 17:00 -- 24:00 & $700$ & $0.063$ & accepted & $0.720$ & accepted\\
	\hline
	& 00:00 -- 02:00 & $86$ & $0.825$ & accepted & $0.709$ & accepted\\
	& 02:00 -- 06:00 & $99$ & $0.451$ & accepted & $0.162$ & accepted\\
	& 06:00 -- 07:00 & $29$ & $ 0.634$ & accepted & $0.538$ & accepted\\
	& 07:00 -- 08:00 & $59$ & $0.424$ & accepted & $0.252$ & accepted\\
	$10^3$ & 08:00 -- 09:00 & $86$ & $0.393$ & accepted & $0.734$ & accepted\\
	& 09:00 -- 15:00 & $811$ & $0.073$ & accepted & $0.955$ & accepted\\
	& 15:00 -- 16:00 & $127$ & $0.928$ & accepted & $0.610$ & accepted\\
	& 16:00 -- 17:00 & $117$ & $0.479$ & accepted & $0.987$ & accepted\\
	& 17:00 -- 18:00 & $111$ & $0.769$ & accepted & $0.516$ & accepted\\
	& 18:00 -- 22:00 & $418$ & $0.660$ & accepted & $0.987$ & accepted\\
	& 22:00 -- 24:00 & $171$ & $0.053$ & accepted & $0.681$ & accepted\\
	\hline
\end{longtable}
\newpage
\begin{longtable}{|c|c|c|c|c|c|c|}
	\caption{Results of the CU KS and dispersion tests (with a significance level of 0.05) applied to each interval of the partition corresponding to the starting point $x^0$. From top to bottom: $m=5,9,17,22,26$. For each interval of each partition, the sample size of the dispersion test is $m$. $H_0$ denotes the null hypothesis of the corresponding test.}\label{tab:Test-initial}\\
	\cline{4-7}
	\multicolumn{3}{c|}{} & \multicolumn{2}{|c|}{CU KS test} & \multicolumn{2}{c|}{Dispersion test} \\
	\cline{1-7}
	$m$ & Interval & $k_i$ & $p$-value & $H_0$ & $p$-value & $H_0$ \\ \hline
	& 00:00 -- 01:00 & $20$ & $0.167$ & accepted & $0.240$ & accepted \\
	& 01:00 -- 02:00 & $11$ & $0.616$ & accepted & $0.151$ & accepted \\
	& 02:00 -- 03:00 & $7$ & $0.887$ & accepted & $0.160$ & accepted \\
	& 03:00 -- 04:00 & $7$ & $0.892$ & accepted & $0.683$ & accepted \\
	& 04:00 -- 05:00 & $12$ & $0.217$ & accepted & $0.856$ & accepted\\
	& 05:00 -- 06:00 & $8$ & $0.426$ & accepted & $0.219$ & accepted\\
	& 06:00 -- 07:00 & $15$ & $0.884$ & accepted & $0.504$ & accepted\\
	& 07:00 -- 08:00 & $27$ & $0.820$ & accepted & $0.164$ & accepted\\
	& 08:00 -- 09:00 & $35$ & $0.875$ & accepted & $0.534$ & accepted\\
	& 09:00 -- 10:00 & $50$ & $0.378$ & accepted & $0.844$ & accepted\\
	& 10:00 -- 11:00 & $48$ & $0.083$ & accepted & $0.884$ & accepted\\ 
	& 11:00 -- 12:00 & $59$ & $0.484$ & accepted & $0.966$ & accepted\\ 
	5& 12:00 -- 13:00 & $51$ & $0.594$ & accepted & $0.765$ & accepted\\ 
	& 13:00 -- 14:00 & $47$ & $0.651$ & accepted & $0.689$ & accepted\\ 
	& 14:00 -- 15:00 & $44$ & $0.817$ & accepted & $0.412$ & accepted\\ 
	& 15:00 -- 16:00 & $45$ & $0.811$ & accepted & $0.168$ & accepted\\ 
	& 16:00 -- 17:00 & $47$ & $0.679$ & accepted & $0.987$ & accepted\\ 
	& 17:00 -- 18:00 & $49$ & $0.486$ & accepted & $0.534$ & accepted\\ 
	& 18:00 -- 19:00 & $37$ & $0.731$ & accepted & $0.344$ & accepted\\ 
	& 19:00 -- 20:00 & $35$ & $0.436$ & accepted & $0.839$ & accepted\\ 
	& 20:00 -- 21:00 & $44$ & $0.904$ & accepted & $0.794$ & accepted\\ 
	& 21:00 -- 22:00 & $43$ & $0.459$ & accepted & $0.693$ & accepted\\ 
	& 22:00 -- 23:00 & $32$ & $0.967$ & accepted & $0.667$ & accepted\\ 
	& 23:00 -- 24:00 & $31$ & $0.306$ & accepted & $0.552$ & accepted\\
	\hline
	& 00:00 -- 01:00 & $33$ & $0.106$ & accepted & $0.527$ & accepted \\
	& 01:00 -- 02:00 & $22$ & $0.658$ & accepted & $0.488$ & accepted \\
	& 02:00 -- 03:00 & $13$ & $0.031$ & rejected & $0.390$ & accepted \\
	& 03:00 -- 04:00 & $14$ & $0.258$ & accepted & $0.857$ & accepted \\
	& 04:00 -- 05:00 & $16$ & $0.441$ & accepted & $0.471$ & accepted\\
	& 05:00 -- 06:00 & $22$ & $0.707$ & accepted & $0.335$ & accepted\\
	& 06:00 -- 07:00 & $23$ & $0.580$ & accepted & $0.608$ & accepted\\
	& 07:00 -- 08:00 & $48$ & $0.500$ & accepted & $0.484$ & accepted\\
	& 08:00 -- 09:00 & $54$ & $0.338$ & accepted & $0.573$ & accepted\\
	& 09:00 -- 10:00 & $97$ & $0.391$ & accepted & $0.886$ & accepted\\
	& 10:00 -- 11:00 & $97$ & $0.149$ & accepted & $0.836$ & accepted\\ 
	& 11:00 -- 12:00 & $108$ & $0.384$ & accepted & $0.999$ & accepted\\ 
	9& 12:00 -- 13:00 & $95$ & $0.911$ & accepted & $0.821$ & accepted\\ 
	& 13:00 -- 14:00 & $82$ & $0.733$ & accepted & $0.923$ & accepted\\ 
	& 14:00 -- 15:00 & $75$ & $0.979$ & accepted & $0.753$ & accepted\\ 
	& 15:00 -- 16:00 & $89$ & $0.909$ & accepted & $0.456$ & accepted\\ 
	& 16:00 -- 17:00 & $82$ & $0.429$ & accepted & $0.923$ & accepted\\ 
	& 17:00 -- 18:00 & $78$ & $0.804$ & accepted & $0.596$ & accepted\\ 
	& 18:00 -- 19:00 & $69$ & $0.277$ & accepted & $0.734$ & accepted\\ 
	& 19:00 -- 20:00 & $69$ & $0.218$ & accepted & $0.477$ & accepted\\ 
	& 20:00 -- 21:00 & $72$ & $0.731$ & accepted & $0.731$ & accepted\\ 
	& 21:00 -- 22:00 & $75$ & $0.449$ & accepted & $0.541$ & accepted\\ 
	& 22:00 -- 23:00 & $60$ & $0.989$ & accepted & $0.681$ & accepted\\ 
	& 23:00 -- 24:00 & $48$ & $0.521$ & accepted & $0.689$ & accepted\\
	\hline
	& 00:00 -- 01:00 & $73$ & $0.729$ & accepted & $0.472$ & accepted \\
	& 01:00 -- 02:00 & $48$ & $0.708$ & accepted & $0.291$ & accepted \\
	& 02:00 -- 03:00 & $39$ & $0.009$ & rejected & $0.010$ & rejected \\
	& 03:00 -- 04:00 & $32$ & $0.203$ & accepted & $0.622$ & accepted \\
	& 04:00 -- 05:00 & $28$ & $0.706$ & accepted & $0.652$ & accepted\\
	& 05:00 -- 06:00 & $38$ & $0.125$ & accepted & $0.607$ & accepted\\
	& 06:00 -- 07:00 & $35$ & $0.908$ & accepted & $0.327$ & accepted\\
	& 07:00 -- 08:00 & $70$ & $0.788$ & accepted & $0.075$ & accepted\\
	& 08:00 -- 09:00 & $121$ & $0.786$ & accepted & $0.577$ & accepted\\
	& 09:00 -- 10:00 & $174$ & $0.421$ & accepted & $0.729$ & accepted\\
	& 10:00 -- 11:00 & $186$ & $0.332$ & accepted & $0.939$ & accepted\\ 
	& 11:00 -- 12:00 & $203$ & $0.474$ & accepted & $0.999$ & accepted\\ 
	17& 12:00 -- 13:00 & $176$ & $0.698$ & accepted & $0.986$ & accepted\\ 
	& 13:00 -- 14:00 & $164$ & $0.589$ & accepted & $0.992$ & accepted\\ 
	& 14:00 -- 15:00 & $161$ & $0.983$ & accepted & $0.570$ & accepted\\ 
	& 15:00 -- 16:00 & $168$ & $0.506$ & accepted & $0.815$ & accepted\\ 
	& 16:00 -- 17:00 & $153$ & $0.361$ & accepted & $0.996$ & accepted\\ 
	& 17:00 -- 18:00 & $149$ & $0.596$ & accepted & $0.528$ & accepted\\ 
	& 18:00 -- 19:00 & $134$ & $0.761$ & accepted & $0.909$ & accepted\\ 
	& 19:00 -- 20:00 & $140$ & $0.101$ & accepted & $0.637$ & accepted\\ 
	& 20:00 -- 21:00 & $141$ & $0.709$ & accepted & $0.760$ & accepted\\ 
	& 21:00 -- 22:00 & $153$ & $0.938$ & accepted & $0.855$ & accepted\\ 
	& 22:00 -- 23:00 & $129$ & $0.887$ & accepted & $0.393$ & accepted\\ 
	& 23:00 -- 24:00 & $94$ & $0.950$ & accepted & $0.296$ & accepted\\
	\hline
	& 00:00 -- 01:00 & $95$ & $0.509$ & accepted & $0.720$ & accepted \\
	& 01:00 -- 02:00 & $70$ & $0.938$ & accepted & $0.529$ & accepted \\
	& 02:00 -- 03:00 & $52$ & $0.008$ & rejected & $0.022$ & rejected \\
	& 03:00 -- 04:00 & $36$ & $0.094$ & accepted & $0.507$ & accepted \\
	& 04:00 -- 05:00 & $34$ & $0.536$ & accepted & $0.420$ & accepted\\
	& 05:00 -- 06:00 & $46$ & $0.045$ & rejected & $0.703$ & accepted\\
	& 06:00 -- 07:00 & $48$ & $0.833$ & accepted & $0.590$ & accepted\\
	& 07:00 -- 08:00 & $83$ & $0.805$ & accepted & $0.062$ & accepted\\
	& 08:00 -- 09:00 & $165$ & $0.576$ & accepted & $0.108$ & accepted\\
	& 09:00 -- 10:00 & $219$ & $0.105$ & accepted & $0.737$ & accepted\\
	& 10:00 -- 11:00 & $235$ & $0.282$ & accepted & $0.960$ & accepted\\ 
	& 11:00 -- 12:00 & $274$ & $0.585$ & accepted & $0.962$ & accepted\\ 
	22& 12:00 -- 13:00 & $233$ & $0.956$ & accepted & $0.984$ & accepted\\ 
	& 13:00 -- 14:00 & $216$ & $0.515$ & accepted & $0.999$ & accepted\\ 
	& 14:00 -- 15:00 & $207$ & $0.872$ & accepted & $0.789$ & accepted\\ 
	& 15:00 -- 16:00 & $213$ & $0.841$ & accepted & $0.905$ & accepted\\ 
	& 16:00 -- 17:00 & $204$ & $0.491$ & accepted & $0.999$ & accepted\\ 
	& 17:00 -- 18:00 & $192$ & $0.534$ & accepted & $0.683$ & accepted\\ 
	& 18:00 -- 19:00 & $173$ & $0.818$ & accepted & $0.968$ & accepted\\ 
	& 19:00 -- 20:00 & $177$ & $0.072$ & accepted & $0.768$ & accepted\\ 
	& 20:00 -- 21:00 & $181$ & $0.655$ & accepted & $0.681$ & accepted\\ 
	& 21:00 -- 22:00 & $196$ & $0.977$ & accepted & $0.810$ & accepted\\ 
	& 22:00 -- 23:00 & $167$ & $0.688$ & accepted & $0.412$ & accepted\\ 
	& 23:00 -- 24:00 & $118$ & $0.963$ & accepted & $0.209$ & accepted\\
	\hline
	& 00:00 -- 01:00 & $112$ & $0.171$ & accepted & $0.679$ & accepted \\
	& 01:00 -- 02:00 & $75$ & $0.933$ & accepted & $0.377$ & accepted \\
	& 02:00 -- 03:00 & $67$ & $0.012$ & rejected & $0.053$ & accepted \\
	& 03:00 -- 04:00 & $46$ & $0.458$ & accepted & $0.450$ & accepted \\
	& 04:00 -- 05:00 & $38$ & $0.987$ & accepted & $0.465$ & accepted\\
	& 05:00 -- 06:00 & $57$ & $0.308$ & accepted & $0.535$ & accepted\\
	& 06:00 -- 07:00 & $56$ & $0.935$ & accepted & $0.739$ & accepted\\
	& 07:00 -- 08:00 & $100$ & $0.882$ & accepted & $0.128$ & accepted\\
	& 08:00 -- 09:00 & $198$ & $0.566$ & accepted & $0.142$ & accepted\\
	& 09:00 -- 10:00 & $259$ & $0.341$ & accepted & $0.844$ & accepted\\
	& 10:00 -- 11:00 & $289$ & $0.091$ & accepted & $0.942$ & accepted\\ 
	& 11:00 -- 12:00 & $320$ & $0.725$ & accepted & $0.984$ & accepted\\ 
	26& 12:00 -- 13:00 & $274$ & $0.915$ & accepted & $0.996$ & accepted\\ 
	& 13:00 -- 14:00 & $257$ & $0.228$ & accepted & $0.999$ & accepted\\ 
	& 14:00 -- 15:00 & $243$ & $0.872$ & accepted & $0.835$ & accepted\\ 
	& 15:00 -- 16:00 & $242$ & $0.574$ & accepted & $0.892$ & accepted\\ 
	& 16:00 -- 17:00 & $236$ & $0.630$ & accepted & $0.942$ & accepted\\ 
	& 17:00 -- 18:00 & $231$ & $0.808$ & accepted & $0.753$ & accepted\\ 
	& 18:00 -- 19:00 & $204$ & $0.682$ & accepted & $0.980$ & accepted\\ 
	& 19:00 -- 20:00 & $209$ & $0.170$ & accepted & $0.830$ & accepted\\ 
	& 20:00 -- 21:00 & $219$ & $0.610$ & accepted & $0.735$ & accepted\\ 
	& 21:00 -- 22:00 & $237$ & $0.803$ & accepted & $0.905$ & accepted\\ 
	& 22:00 -- 23:00 & $198$ & $0.614$ & accepted & $0.366$ & accepted\\ 
	& 23:00 -- 24:00 & $147$ & $0.972$ & accepted & $0.032$ & accepted\\
	\hline
\end{longtable}
\newpage
\begin{longtable}{|c|c|c|c|c|c|c|}
	\caption{Results of the CU KS and dispersion tests (with a significance level of 0.05) applied to each interval of the final (infeasible) partition obtained by solving problem \ref{optBestFit} for different values of the parameter $m$, with $w$ fixed to 1. From top to bottom: $m=5,9,17,22,26$. For each interval of each partition, the sample size of the dispersion test is $m$. $H_0$ denotes the null hypothesis of the corresponding test.}\label{tab:Test-infeasible-final}\\
	\cline{4-7}
	\multicolumn{3}{c|}{} & \multicolumn{2}{|c|}{CU KS test} & \multicolumn{2}{c|}{Dispersion test} \\
	\cline{1-7}
	$m$ & Interval & $k_i$ & $p$-value & $H_0$ & $p$-value & $H_0$ \\ \hline
	& 00:00 -- 06:00 & $65$ & $0.068$ & accepted & $0.472$ & accepted \\
	& 06:00 -- 07:00 & $15$ & $0.884$ & accepted & $0.504$ & accepted \\
	& 07:00 -- 08:00 & $27$ & $0.820$ & accepted & $0.164$ & accepted \\
	& 08:00 -- 09:00 & $35$ & $0.875$ & accepted & $0.534$ & accepted\\
	& 09:00 -- 10:00 & $50$ & $0.378$ & accepted & $0.844$ & accepted\\
	5 & 10:00 -- 12:00 & $107$ & $0.734$ & accepted & $0.938$ & accepted\\
	& 12:00 -- 13:00 & $51$ & $0.594$ & accepted & $0.765$ & accepted\\
	& 13:00 -- 14:00 & $47$ & $0.651$ & accepted & $0.689$ & accepted\\
	& 14:00 -- 15:00 & $44$ & $0.817$ & accepted & $0.412$ & accepted\\
	& 15:00 -- 24:00 & $363$ & $0.214$ & accepted & $0.568$ & accepted\\
	\hline
	& 00:00 -- 02:00 & $55$ & $0.249$ & accepted & $0.607$ & accepted \\
	& 02:00 -- 04:00 & $27$ & $0.309$ & accepted & $0.501$ & accepted \\
	& 04:00 -- 05:00 & $16$ & $0.441$ & accepted & $0.471$ & accepted\\
	& 05:00 -- 06:00 & $22$ & $0.707$ & accepted & $0.335$ & accepted \\
	& 06:00 -- 07:00 & $23$ & $0.580$ & accepted & $0.608$ & accepted\\
	& 07:00 -- 08:00 & $48$ & $0.500$ & accepted & $0.484$ & accepted\\
	9 & 08:00 -- 09:00 & $54$ & $0.338$ & accepted & $0.573$ & accepted\\
	& 09:00 -- 16:00 & $643$ & $0.060$ & accepted & $0.717$ & accepted\\
	& 16:00 -- 17:00 & $82$ & $0.429$ & accepted & $0.923$ & accepted\\
	& 17:00 -- 18:00 & $78$ & $0.804$ & accepted & $0.596$ & accepted\\
	& 18:00 -- 22:00 & $285$ & $0.919$ & accepted & $0.989$ & accepted\\
	& 22:00 -- 23:00 & $60$ & $0.989$ & accepted & $0.681$ & accepted\\
	& 23:00 -- 24:00 & $48$ & $0.522$ & accepted & $0.689$ & accepted\\
	\hline
	& 00:00 -- 02:00 & $121$ & $0.094$ & accepted & $0.535$ & accepted \\
	& 02:00 -- 05:00 & $99$ & $0.098$ & accepted & $0.067$ & accepted \\
	& 05:00 -- 07:00 & $73$ & $0.650$ & accepted & $0.203$ & accepted\\
	& 07:00 -- 08:00 & $70$ & $0.788$ & accepted & $0.075$ & accepted \\
	& 08:00 -- 09:00 & $121$ & $0.786$ & accepted & $0.577$ & accepted\\
	& 09:00 -- 10:00 & $174$ & $0.421$ & accepted & $0.729$ & accepted\\
	17 & 10:00 -- 14:00 & $729$ & $0.089$ & accepted & $0.995$ & accepted\\
	& 14:00 -- 16:00 & $329$ & $0.982$ & accepted & $0.410$ & accepted\\
	& 16:00 -- 17:00 & $153$ & $0.361$ & accepted & $0.996$ & accepted\\
	& 17:00 -- 18:00 & $149$ & $0.596$ & accepted & $0.528$ & accepted\\
	& 18:00 -- 22:00 & $568$ & $0.586$ & accepted & $0.926$ & accepted\\
	& 22:00 -- 24:00 & $223$ & $0.071$ & accepted & $0.793$ & accepted\\
	\hline
	& 00:00 -- 02:00 & $165$ & $0.198$ & accepted & $0.743$ & accepted\\
	& 02:00 -- 06:00 & $168$ & $0.117$ & accepted & $0.122$ & accepted\\
	& 06:00 -- 07:00 & $48$ & $0.833$ & accepted & $0.590$ & accepted\\
	& 07:00 -- 08:00 & $83$ & $0.805$ & accepted & $0.062$ & accepted\\
	& 08:00 -- 09:00 & $165$ & $0.576$ & accepted & $0.108$ & accepted\\
	& 09:00 -- 10:00 & $219$ & $0.105$ & accepted & $0.737$ & accepted\\
	22& 10:00 -- 14:00 & $958$ & $0.097$ & accepted & $0.994$ & accepted\\
	& 14:00 -- 16:00 & $420$ & $0.952$ & accepted & $0.561$ & accepted\\
	& 16:00 -- 17:00 & $204$ & $0.491$ & accepted & $0.999$ & accepted\\
	& 17:00 -- 18:00 & $192$ & $0.534$ & accepted & $0.683$ & accepted\\
	& 18:00 -- 22:00 & $772$ & $0.436$ & accepted & $0.968$ & accepted\\
	& 22:00 -- 23:00 & $167$ & $0.688$ & accepted & $0.412$ & accepted\\
	& 23:00 -- 24:00 & $118$ & $0.963$ & accepted & $0.209$ & accepted\\
	\hline
	& 00:00 -- 01:00 & $112$ & $0.171$ & accepted & $0.679$ & accepted\\
	& 01:00 -- 02:00 & $75$ & $0.933$ & accepted  & $0.378$ & accepted\\
	& 02:00 -- 06:00 & $208$ & $0.072$ & accepted & $0.080$ & accepted\\
	& 06:00 -- 07:00 & $56$ & $0.935$ & accepted & $0.739$ & accepted\\
	& 07:00 -- 08:00 & $100$ & $0.882$ & accepted & $0.128$ & accepted\\
	& 08:00 -- 09:00 & $198$ & $0.566$ & accepted & $0.142$ & accepted\\
	& 09:00 -- 10:00 & $259$ & $0.341$ & accepted & $0.844$ & accepted\\
	& 10:00 -- 11:00 & $289$ & $0.091$ & accepted & $0.942$ & accepted\\
	26& 11:00 -- 12:00 & $320$ & $0.725$ & accepted & $0.984$ & accepted\\
	& 12:00 -- 13:00 & $274$ & $0.915$ & accepted & $0.996$ & accepted\\
	& 13:00 -- 15:00 & $500$ & $0.439$ & accepted & $0.971$ & accepted\\
	& 15:00 -- 16:00 & $242$ & $0.574$ & accepted & $0.892$ & accepted\\
	& 16:00 -- 18:00 & $467$ & $0.895$ & accepted & $0.939$ & accepted\\
	& 18:00 -- 21:00 & $632$ & $0.643$ & accepted & $0.950$ & accepted\\
	& 21:00 -- 22:00 & $237$ & $0.803$ & accepted & $0.905$ & accepted\\
	& 22:00 -- 24:00 & $345$ & $0.034$ & rejected & $0.440$ & accepted\\
	\hline
\end{longtable}

%
%

%
\section*{Conflict of interest}

The authors declare that they have no conflict of interest.



\end{document}